\documentclass{amsart}
\usepackage{amssymb,amsmath,amsfonts}

\newtheorem{theorem}{Theorem}[section] 
\newtheorem{definition}[theorem]{Definition}
\newtheorem{proposition}[theorem]{Proposition}
\newtheorem{lemma}[theorem]{Lemma}
\newtheorem{corollary}[theorem]{Corollary}

\theoremstyle{definition}
\newtheorem{remark}[theorem]{Remark}
\newtheorem{example}[theorem]{Example}

\def\A{\mathbb{A}}
\def\CC{\mathbb{C}}
\def\F{\mathcal{F}}
\def\FF{\mathbb{F}}
\def\G{\mathbb{G}}
\def\H{\mathcal{H}}
\def\M{\mathbb{M}}
\def\P{\mathbb{P}}
\def\Q{\mathbb{Q}}
\def\R{\mathcal{R}}
\def\RR{\mathbb{R}}
\def\U{\mathbb{U}}
\def\Z{\mathbb{Z}}

\def\<{\ensuremath{\langle}}
\def\>{\ensuremath{\rangle}}

\DeclareMathOperator{\Hom}{Hom}
\DeclareMathOperator{\init}{in}
\DeclareMathOperator{\PGL}{PGL}
\DeclareMathOperator{\Spec}{Spec}
\DeclareMathOperator{\Trop}{Trop} 

\begin{document}

\title{Realization spaces for tropical fans}

\author[Katz]{Eric Katz}
\address{Department of Mathematics, University of Texas at Austin, Austin, TX 78712}
\email{eekatz@math.utexas.edu}
\author[Payne]{Sam Payne$^*$}\thanks{$^*$Supported by the Clay Mathematics Institute.}
\address{Stanford University, Mathematics, Bldg. 380, 450 Serra Mall, Stanford, CA 94305}
\email{spayne@stanford.edu}

\begin{abstract}
We introduce a moduli functor for varieties whose tropicalization realizes a given weighted fan and show that this functor is an algebraic space in general, and is represented by a scheme when the associated toric variety is quasiprojective.  We study the geometry of these tropical realization spaces for the matroid fans studied by Ardila and Klivans, and show that the tropical realization space of a matroid fan is a torus torsor over the realization space of the matroid.  As a consequence, we deduce that these tropical realization spaces satisfy Murphy's Law.
\end{abstract}

\maketitle


\section{Introduction}

The tropicalization of a $d$-dimensional subvariety of a torus with respect to the trivial valuation is the underlying set of a pure $d$-dimensional fan with a locally constant positive integer weight function on its smooth points.  The question of which subvarieties of a torus, if any, have a given tropicalization has generated many discussions, but the literature treats only the special cases of existence of realizations for curves of genus zero or one \cite{Speyer07b} and tropical linear spaces \cite{Mikhalkin08}, in characteristic zero.  Here we use Hilbert schemes and Tevelev's approach to compactifications of subvarieties of tori \cite{Tevelev07} to address the general case, constructing a fine moduli scheme over the integers that parametrizes varieties whose tropicalization realizes a given weighted simplicial fan.

Let $T$ be the torus whose lattice of one parameter subgroups is $N \cong \Z^n$, and let $N_\RR = N \otimes_\Z \RR$.  Fix a pure $d$-dimensional simplicial fan $\Delta$ in $N_\RR$ with a weight function $w$ that assigns a positive integer $w(\sigma)$ to each $d$-dimensional cone $\sigma$ in $\Delta$.  Say that a \emph{tropical realization} of $(\Delta,w)$ is a subvariety $Y^\circ$ of $T$ whose tropicalization is the underlying set $|\Delta|$, with weight function given by $w(\sigma)$ on the relative interior of $\sigma$.   

In Section~\ref{compactifications}, we show that a subvariety $Y^\circ$ of $T$ is a tropical realization of $(\Delta, w)$ if and only if its closure $Y$ in the toric variety $X(\Delta)$ is proper and the intersection number $(V(\sigma) \cdot Y)$ is equal to $w(\sigma)$ for every $d$-dimensional cone $\sigma$ in $\Delta$.  The existence of Hilbert schemes makes it convenient to work with these proper subvarieties of toric varieties that appear as closures of tropical realizations rather than with the tropical realizations themselves, especially when working in families.

\begin{definition}
Let $S$ be a scheme.  A \emph{family of tropical realizations} of $(\Delta,w)$ over $S$ is a subscheme $Y$ of $X(\Delta) \times S$, flat and proper over $S$, such that every geometric fiber is reduced, irreducible, and has intersection number $w(\sigma)$ with $V(\sigma)$, for every $d$-dimensional cone $\sigma$ in $\Delta$.
\end{definition}

\noindent The pullback of a family of tropical realizations of $\Delta$ under a morphism $S' \rightarrow S$ is a family of tropical realizations over $S'$, so the map $\R_{\Delta,w}$ taking a scheme $S$ to the set of all families of tropical realizations of $\Delta$ over $S$ is a contravariant functor from schemes to sets.  In Section~\ref{representability} we show that this functor is an algebraic space, in the sense of \cite{Knutson71}.  Recall that the toric variety $X(\Delta)$ is quasiprojective if and only if the fan $\Delta$ can be extended to a regular subdivision of $N_\RR$.

\begin{theorem} \label{representable}
If the toric variety $X(\Delta)$ is quasiprojective then the tropical realization functor $\R_{\Delta,w}$ is represented by a scheme of finite type.
\end{theorem}

\noindent In other words, there is a scheme $\mathbf{R}_{\Delta,w}$ of finite type and a universal family of tropical realizations of $(\Delta,w)$ over $\mathbf{R}_{\Delta, w}$ such that any family of tropical realizations of $\Delta$ over an arbitrary scheme $S$ is the pullback of this universal family under a unique morphism $S \rightarrow \mathbf{R}_{\Delta, w}$.  We call $\mathbf{R}_{\Delta,w}$ the \emph{tropical realization space} of $(\Delta, w)$.  Theorem~\ref{representable} has been applied by Katz to construct nonarchimedean analytic moduli spaces for realizations of tropical polyhedral complexes over valued fields \cite{Katz10}.

\begin{remark}
If $(\Delta', w')$ is a refinement of $(\Delta, w)$ then a subvariety of $T$ is a tropical realization of $(\Delta', w')$ if and only if it is a tropical realization of $(\Delta, w)$.  So there is a natural bijection between the $k$-points of $\mathbf R_{\Delta, w}$ and of $\mathbf R_{\Delta', w'}$, for every field $k$, although these schemes need not be isomorphic.
\end{remark}

In Section~\ref{matroids}, we study the tropical realization spaces of a special class of fans associated to matroids.  The \emph{matroid fan} $\Delta_\M$ associated to a matroid $\M$ was introduced by Ardila and Klivans \cite{ArdilaKlivans06}, who called it the fine subdivision of the Bergman fan of the matroid.

\begin{theorem} \label{smooth and surjective}
Let $\M$ be a matroid.  The tropical realization space of the matroid fan $\Delta_\M$ is naturally a torus torsor of rank $\# \M-c(\M)$ over the realization space of $\M$.
\end{theorem}

\noindent Here $\# \M$ is the number of elements in the underlying set of the matroid and $c(\M)$ is the number of connected components of $\M$.

\begin{corollary} \label{realizability}
The matroid fan $\Delta_\M$ is tropically realizable over a field $k$ if and only if the matroid $\M$ is realizable over $k$.  
\end{corollary}

\noindent For instance, the matroid associated to the configuration of all seven lines in the projective plane $\P^2(\mathbb F_2)$ over the field with two elements is realizable only over fields of characteristic two.  The associated matroid fan is a two-dimensional fan in a six-dimensional real vector space with twenty-one maximal cones corresponding to the twenty-one full flags in $\P^2(\mathbb F_2)$, and it is the tropicalization of a two-dimensional subvariety of a six-dimensional torus only over fields of characteristic two.    

\begin{remark}
Corollary~\ref{realizability} has been known as a folklore result in the tropical geometry community, and a version of this result over the complex numbers has been reported by Mikhalkin, based on joint work with Sturmfels and Ziegler \cite{Mikhalkin08}, but there is no proof in the literature.  The geometric tropicalization theory of Hacking, Keel, and Tevelev \cite{HackingKeelTevelev09} implies that $\Delta_\M$ is the tropicalization of the complement of a hyperplane arrangement over a field $k$ if and only if that hyperplane arrangement realizes the matroid $\M$.  In Section~\ref{matroids}, we show that every tropical realization of $\Delta_\M$ is the complement of a hyperplane arrangement, which proves the corollary.  For those interested only in the corollary, this is simpler than deducing it from Theorem~\ref{smooth and surjective}, whose proof is more technical.
\end{remark}

Combining Theorem~\ref{smooth and surjective} and Mn\"ev's Universality Theorem, we deduce that tropical realization spaces can have arbitrary singularity types, even when the tropical realizations themselves are all smooth.  Following Vakil \cite{Vakil06}, we say that a singularity type is an equivalence class of pointed schemes $(Y,y)$, with the equivalence relation generated by setting $(Y,y) \sim (Y',y')$ if there is a smooth morphism from $Y$ to $Y'$ that takes $y$ to $y'$.  A collection of schemes satisfies Murphy's Law if every singularity type occurs on some scheme in the collection.

We say that a weighted fan $(\Delta,w)$ is \emph{multiplicity free} if $w(\sigma)$ is equal to one for all $\sigma$.  In this case, we write $\mathbf{R}_\Delta$ for the tropical realization space $\mathbf{R}_{\Delta,w}$.

\begin{theorem} \label{Murphy}
Every singularity type of finite type over $\Spec \Z$ appears in the tropical realization space of a multiplicity free two-dimensional fan $\Delta$ such that the universal family is smooth over $\mathbf{R}_\Delta$.
\end{theorem}

\noindent In other words, tropical realization spaces of multiplicity free two-dimensional fans satisfy Murphy's Law.

Mn\"ev's Universality Theorem also implies that deciding whether an arbitrary matroid is realizable over a given field $k$ that is not algebraically closed is as difficult as the generalized Hilbert's Tenth Problem of deciding whether an arbitrary system of polynomials with integer coefficients has a solution over $k$ \cite{Sturmfels87}.  This problem is known to be undecidable over many function fields, such as  $\FF_q(t)$, $\RR(t)$, and $\CC(t,u)$; decidability over the rational numbers is an open problem \cite{Poonen08}.

\smallskip

\noindent \textbf{Conventions.}  Throughout this paper, all tropicalizations are with respect to the trivial valuation and all schemes are locally of finite type.

\smallskip

\noindent \textbf{Acknowledgments.}  We are grateful to D. Eisenbud, E. Feichtner, D. Helm, and R. Vakil for helpful discussions related to this work and to the referee for several useful comments and suggestions.

\section{Realizations and compactifications} \label{compactifications}

As in the introduction, $\Delta$ is a pure $d$-dimensional fan in $N_\RR$ and $X = X(\Delta)$ is the associated toric variety with dense torus $T$.  The weight function $w$ assigns a positive integer $w(\sigma)$ to each $d$-dimensional cone $\sigma$ in $\Delta$.  Here we discuss how subvarieties of $T$ that are tropical realizations of $(\Delta, w)$ can be characterized in terms of their closures in $X$.

Let $Y^\circ$ be a $d$-dimensional subvariety of $T$.  Recall that the tropicalization $\Trop(Y^\circ)$ is the set of vectors $v \in N_\RR$ such that the initial degeneration $\init_v(Y^\circ)$ in $T$ is nonempty.  It is the underlying set of a rational fan of pure dimension $d$, and comes with a tropical multiplicity function $m$ on its subset of smooth points, defined as follows.  We say that a point in $\Trop(Y^\circ)$ is \emph{smooth} if it has an open neighborhood homeomorphic to $\RR^d$.   There is an open dense set of smooth points $v$ in $\Trop(Y^\circ)$ such that a $d$-dimensional subtorus $T_v$ acts freely on $\init_v Y^\circ$.  For such $v$, the tropical multiplicity $m(v)$ is defined to be the length of the zero dimensional quotient $\init_v Y^\circ/ T_v$.  The balancing condition implies that the tropical multiplicity extends to a locally constant function on the set of all smooth points of $\Trop(Y^\circ)$.  In particular, the tropical multiplicity is constant on the relative interior of each maximal cone of $\Delta$.  The balancing condition then says that these integers associated to maximal cones are a Minkowski weight, in the sense of \cite{FultonSturmfels97}.  Since the balancing condition is necessary for a weighted fan to be tropically realizable, nothing is lost by assuming that all weighted fans in this paper are balanced.

We now state the conventional definition of a tropical realization over a field $k$.  In Proposition~\ref{intersection criterion} below, we show that the tropical realizations in this naive sense are exactly those subvarieties of $T$ whose closure in $X(\Delta)$ is a family of tropical realizations over $\Spec k$.  To avoid any possibility of confusion, we always write $Y^\circ$ for a subvariety of the torus $T$, and $Y$ for its closure in $X(\Delta)$.

\begin{definition}
A \emph{tropical realization} of $(\Delta,w)$ is a subvariety $Y^\circ$ of $T$ such that the underlying set of $\Trop(Y^\circ)$ is $|\Delta|$ and the tropical multiplicity $m(v)$ is equal to $w(\sigma)$ for $v$ in the relative interior of $\sigma$.
\end{definition}

\begin{lemma}
Let $Y^\circ$ be a $d$-dimensional subvariety of $T$.  The underlying set of $\Trop(Y^\circ)$ is equal to $|\Delta|$ if and only if its closure $Y$ in $X(\Delta)$ is proper and meets every $T$-invariant subvariety.
\end{lemma}

\begin{proof}
This follows from Lemma~2.2 of \cite{Tevelev07} in the case where $X$ is smooth, since $Y$ and $|\Delta|$ both have pure dimension $d$.  For the general case, and for further references, see Section~8 of \cite{Katz09}.
\end{proof}

Since $\Delta$ is simplicial, the toric variety $X(\Delta)$ is the coarse moduli space of a smooth toric stack \cite{Iwanari09}. It follows that there is a natural intersection product on the rational Chow groups $A_*(X)_\Q$ \cite{Vistoli89} with all of the refined properties of the intersection product on Chow groups of smooth varieties, as presented in Chapter~8 of \cite{IT}.  In particular, if $\sigma$ is a $d$-dimensional cone in $\Delta$, then $V(\sigma) \cdot Y$ is a well-defined class in $A_0(V(\sigma) \cap Y)_\Q$.  If $Y$ is proper then this class has a well-defined degree, which is a rational number denoted $(V (\sigma) \cdot Y)$. 

\begin{lemma}
If $\Trop(Y^\circ) = |\Delta|$ then the intersection number $(V(\sigma) \cdot Y)$ is equal to the tropical multiplicity $m(v)$ for $v$ in the relative interior of $\sigma$.
\end{lemma}

\begin{proof}
If $\Delta$ is sufficiently fine, then $Y$ is Cohen-Macaulay and the scheme theoretic intersection $V(\sigma) \cap Y$ is canonically identified with the zero-dimensional quotient $\init_v(Y^\circ)/T_v$, for a generic $v$ in the relative interior of $\sigma$.  The general case follows by choosing a suitable refinement and applying the projection formula.  See Section~9 of \cite{Katz09} for details.
\end{proof}

\noindent  We now prove the following proposition, which was mentioned in the introduction.

\begin{proposition} \label{intersection criterion}
A subvariety $Y^\circ$ of the torus $T$ is a tropical realization of $\Delta$ if and only if its closure $Y$ in $X(\Delta)$ is proper, and
\[
(V(\sigma) \cdot Y) = w(\sigma),
\]
for every $d$-dimensional cone $\sigma \in \Delta$.  
\end{proposition}

\begin{proof} 
If $Y^\circ$ is a tropical realization of $(\Delta,w)$ then $Y$ is proper and $(V(\sigma) \cdot Y) = w(\sigma)$ by the two preceding lemmas.  We now show the converse.  

Suppose $Y$ is proper and $(V(\sigma) \cdot Y) = w(\sigma)$ for every $d$-dimensional cone $\sigma$ in $\Delta$.  Since $Y$ is proper, and both $Y$ and $\Delta$ are $d$-dimensional, $\Trop(Y^\circ)$ is a union of $d$-dimensional cones of $\Delta$.  For any $d$-dimensional cone $\sigma$ in $\Delta$, the intersection number $(V(\sigma) \cdot Y)$ is positive, so $Y$ must intersect $V(\sigma)$.  It follows that $\Trop(Y^\circ)$ contains $\sigma$, and since this holds for all maximal cones of $\Delta$, $\Trop(Y^\circ) = |\Delta|$.  By the preceding lemma, $m(v)$ is equal to $(V(\sigma) \cdot Y)$ for $v$ in the relative interior of $\sigma$, and the proposition follows.
\end{proof}

\noindent The proposition motivates our definitions of a family of tropical realizations of $(\Delta, w)$ and the tropical realization functor $\R_{\Delta,w}$ in the introduction.

\begin{corollary}
Let $k$ be a field.  Then $\R_{\Delta,w}(k)$ is the set of closures in $X(\Delta)$ of the tropical realizations of $(\Delta,w)$ over $k$.
\end{corollary}

\begin{proof}
The corollary follows from Proposition~\ref{intersection criterion}, since $\R_{\Delta,w}(k)$ is the set of proper $k$-subvarieties of $X(\Delta)$ whose intersection number with $V(\sigma)$ is $w(\sigma)$ for every $d$-dimensional cone $\sigma$ in $\Delta$.\end{proof}

\section{Representability} \label{representability}

In Section~\ref{compactifications} we showed that a subvariety $Y^\circ$ of $T$ is a tropical realization of $(\Delta,w)$ if and only if its closure $Y$ in $X(\Delta)$ is proper and $(V(\sigma) \cdot Y) = w(\sigma)$ for every $d$-dimensional cone $\sigma$ in $\Delta$.  

\begin{proposition} \label{algebraic space}
The tropical realization functor $\R_{\Delta,w}$ is an open algebraic subspace of the Hilbert functor of any toric compactification of $X(\Delta)$.
\end{proposition}

\begin{proof}
Let $\Sigma$ be a complete simplicial fan that contains $\Delta$ as a subfan; for a construction of such a fan, see Theorem~III.2.8 of \cite{Ewald96}.  Then $X(\Sigma)$ is a proper simplicial toric variety that contains $X(\Delta)$ as a $T$-invariant open subvariety.  The Hilbert functor of $X(\Sigma)$ is an algebraic space \cite{Artin69}, and there is an open subfunctor of the Hilbert functor that parametrizes reduced and irreducible $d$-dimensional subschemes that are contained in $X(\Delta)$.  Since intersection numbers are locally constant in flat families, it follows that $\R_{\Delta,w}$ is also an open subfunctor of the Hilbert functor of $X(\Sigma)$, and hence is an algebraic space.
\end{proof}

We now prove Theorem~\ref{representable}, which says that the tropical realization functor $\R_{\Delta,w}$ is represented by a scheme of finite type over the integers when $X(\Delta)$ is quasiprojective.  Recall that toric varieties are canonically defined over the integers, and $X(\Delta)$ is quasiprojective if and only if $\Delta$ can be extended to a regular subdivision of $N_\RR$, which means that $X(\Delta)$ can be embedded as a $T$-invariant open subvariety of a projective toric variety. 

\begin{theorem} \label{finite type}
If $X(\Delta)$ is quasiprojective then the tropical realization space $\mathbf{R}_{\Delta,w}$ is represented by a scheme of finite type.
\end{theorem}

\begin{proof}
Fix a projective toric variety $X(\Sigma) \subset \P^r$ that contains $X(\Delta)$ as a $T$-invariant open subvariety.  By Proposition~\ref{algebraic space}, the tropical realization functor $\R_{\Delta,w}$ is represented by an open subscheme $\mathbf R_{\Delta,w}$ of the Hilbert scheme of $X(\Sigma)$.  We must show that $\mathbf R_{\Delta,w}$ is of finite type.  Suppose $Y$ is the closure of a tropical realization of $(\Delta,w)$.  We will show that there are only finitely many possibilities for the Hilbert polynomial of $Y$, and the theorem follows because the Hilbert scheme parametrizing subschemes of $X(\Sigma)$ with fixed Hilbert polynomial is of finite type.

First, we claim that the degree of $Y$ in $\P^r$ is determined by $w$.  Now the Chow homology class $[Y]$ in $A_d(X(\Sigma))_\Q$ is determined by the intersection numbers $(V(\sigma) \cdot Y)$ for $d$-dimensional cones $\sigma$ in $\Sigma$ \cite{FultonSturmfels97}, and $Y$ is disjoint from $V(\sigma)$ if $\sigma$ is not contained in $\Delta$.  Therefore $(V(\sigma) \cdot Y)$ is $w(\sigma)$ if $\sigma$ is in $\Delta$, and zero otherwise.  This determines $[Y]$ and the claim follows, since the degree of $Y$ is given by the push forward of $[Y]$ under the given embedding of $X(\Sigma)$ in $\P^r$.

Next, the Castelnuovo-Mumford regularity of $Y$ is bounded in terms of its degree and $r$ \cite{BayerMumford93}.  So there is an integer $N_0$, depending only on $(\Delta,w)$, such that the Hilbert function of $Y$ evaluated at $N$ agrees with the Hilbert polynomial for $N \geq N_0$.  Therefore, the Hilbert polynomial of $Y$ is determined by the value of the Hilbert function at $N_0, \ldots, N_0 + d$.

Finally, the Hilbert function of a $d$-dimensional subvariety of $\P^r$ is bounded uniformly in terms of $d$, $r$, and its degree; for the construction of such a bound, see Exercise 3.28 in \cite{Kollar96}.  It follows that there are only finitely many possible values for the Hilbert function of $Y$ at $N_0, \ldots, N_0 + d$, and hence there are only finitely many possibilities for the Hilbert polynomial of $Y$, as required.
\end{proof}

\noindent Since any fan $\Delta$ has a refinement $\Delta'$ that can be extended to a regular subdivision of $N_\RR$, Theorem~\ref{finite type} implies that, for an arbitrary weighted simplicial fan $(\Delta,w)$, the tropical realizations of $(\Delta,w)$ over a field $k$ are in natural bijection with the $k$-points of a scheme of finite type.

\begin{remark}
These tropical realization spaces can also be used to construct more general moduli, such as for tropical realizations of pairs.  If $\Delta'$ is a subfan of positive codimension in $\Delta$, then inside the product
\[
\R_{\Delta',w'} \times \R_{\Delta,w}
\]
there is the closed algebraic subspace parametrizing pairs $(Y',Y)$ such that $Y'$ is contained in $Y$.  The fibers of this moduli space of pairs under second projection parametrize tropical realizations of $(\Delta',w')$ that are contained in a fixed subvariety $Y$ of $X(\Delta)$.
\end{remark}

\section{Tropical realizations of matroid fans} \label{matroids}

Let $\M$ be a loop free matroid on a finite set $E$, and let $e_0, \ldots, e_n$ be the elements of $E$.  Let $N$ be the lattice
\[
N = \Z^E / \< e_0 + \cdots + e_n \>.
\]
The \emph{matroid fan} $\Delta_\M$ is a simplicial fan in $N_\RR$ that encodes the lattice of flats of $\M$.  Ardila and Klivans introduced this fan in \cite{ArdilaKlivans06} and called it the fine subdivision of the Bergman fan of the matroid; the fan is defined as follows.  For a subset $I \subset E$, let $e_I$ be the vector 
\[
e_I = \sum_{e_i \in I} e_i
\]
in $N_\RR$.

The rays of the matroid fan $\Delta_\M$ correspond to the proper flats $F \subsetneq E$ of the matroid, and the ray $\rho_F$ corresponding to a flat $F$ is spanned by $e_F$.  More generally, the $k$-dimensional cones of the matroid fan correspond to the $k$-step flags of proper flats.  If $\F$ is a flag of flats $F_1 \subset \cdots \subset F_k$ then the corresponding cone $\sigma_\F$ is spanned by $\{ e_{F_1}, \ldots, e_{F_k} \}$.  Since each cone $\sigma_\F$ in $\Delta_\M$ is spanned by a subset of a basis for the lattice $N$, the toric variety $X(\Delta_\M)$ is smooth.  Furthermore, since every flag of flats in a matroid can be extended to a full flag, the matroid fan $\Delta_\M$ is of pure dimension equal to the rank of $\M$ minus one.

\begin{example} \label{uniform matroid}
Let $\U_n$ be the uniform matroid on $\{0, \ldots, n\}$, the matroid in which every subset is a flat.  Then the matroid fan $\Delta_{\U_n}$ in $\RR^{n+1}/(1, \ldots, 1)$ is the first barycentric subdivision of the fan corresponding to $\P^n$, and $X(\Delta_{\U_n})$ is obtained from $\P^n$ by a sequence of blowups
\[
X(\Delta_{\U_n}) = X_{n-1} \rightarrow \cdots \rightarrow X_1 \rightarrow X_0 = \P^n,
\]
where $X_{i+1} \rightarrow X_i$ is the blowup along the strict transforms of the $i$-dimensional $T$-invariant subvarieties of $\P^n$.  The matroid fan $\Delta_{\U_n}$ can also be realized as the normal fan of the $n$-dimensional permutahedron.
\end{example}

\noindent Note that the labeling $E = \{e_0, \ldots, e_n\}$ of the underlying set of the matroid $\M$ induces an inclusion of the matroid fan $\Delta_\M$ as a subfan of $\Delta_{\U_n}$.  In particular, the toric variety $X(\Delta_\M)$ is quasiprojective, so the tropical realization functor $\R_{\Delta_\M}$ is represented by a scheme of finite type.  Furthermore, the dense torus $T$ in $X(\Delta_\M)$ is naturally identified with the dense torus in $\P^n$.

\begin{proposition} \label{linear subspace}
Let $Y^\circ$ be a tropical realization of $\Delta_\M$.  Then the closure of $Y^\circ$ in $\P^n$ is a $d$-dimensional linear subspace.
\end{proposition}

We will prove the proposition using a projection to $\P^d$ and the following lemmas. After possibly renumbering the elements of $E$, we may assume $\{e_0, \ldots, e_d\}$ is a basis for $\M$.  Let $\U$ be the uniform matroid $\U_d$.  Consider the projection from $N_\RR$ to $\RR^{d+1}/(1, \ldots, 1)$ taking $e_i$ to the image of the $i$th standard basis vector for $i \leq d$, and taking $e_j$ to zero for $j > d$.  Every cone in $\Delta_\M$ projects into some cone of $\Delta_\U$, inducing a map of fans $\Delta_\M \rightarrow \Delta_\U$.  Composing the induced map of toric varieties $X(\Delta_\M) \rightarrow X(\Delta_\U)$  with the birational projection $X(\Delta_\U) \rightarrow \P^d$ described in Example~\ref{uniform matroid} gives a natural map
\[
p: X(\Delta_\M) \rightarrow \P^d.
\]
This map can also be factored as a birational morphism followed by a linear projection, as follows.  Let $D_i$ be the $i$th coordinate hyperplane in $\P^n$.  Then $p$ factors as
\[
X(\Delta_\M) \xrightarrow{p_1} \P^n \smallsetminus (D_0 \cap \cdots \cap D_d) \xrightarrow{p_2} \P^d,
\]
where $p_1$ is the birational morphism induced by the open immersion of $X(\Delta_\M)$ in the iterated blowup $X(\Delta_{U_n})$ of $\P^n$, and $p_2$ is the linear projection away from $D_0 \cap \cdots \cap D_d$.

Let $Y$ be the closure in $X(\Delta_\M)$ of a tropical realization $Y^\circ$ of $\Delta_\M$.

\begin{lemma}
The map $p$ takes $Y$ birationally onto $\P^d$.
\end{lemma}

\begin{proof}
We show that the map from $\Delta_\M$ to $\Delta_\U$ is of tropical degree one, and it follows that the map $Y \rightarrow X(\Delta_\U)$ is birational, by Theorem~1.1 of \cite{SturmfelsTevelev08}.  The lemma follows by composing this map with the birational projection from $X(\Delta_\U)$ to $\P^d$.

Let $v$ be a point in the interior of the maximal cone $\sigma$ in $\Delta_\U$ corresponding to a full flag $S_1 \subset \cdots \subset S_d$ of proper subsets of $\{0, \ldots, d\}$.  Let $\F$ be the flag $F_1 \subset \cdots \subset F_d$ in which $F_i$ is the flat spanned by $\{e_j \ | \ j \in S_i\}$.  Then $\sigma_\F$ is the unique cone in $\Delta_\M$ whose image meets the relative interior of $\sigma$, and it projects bijectively onto $\sigma$.  Therefore, the preimage of $v$ in $|\Delta_\M|$ is a single point in the relative interior of $\sigma_\F$.  Since the lattice $N \cap \mathrm{span}(\sigma_\F)$ maps isomorphically onto $\Z^{d+1}/(1, \ldots, 1)$, the map of fans $\Delta_\M \rightarrow \Delta_\U$ is of tropical degree one, as required.
\end{proof}

\noindent Since $p$ is birational and $\P^d$ is normal, the complement of the open set over which $p$ is an isomorphism has codimension at least two.  In particular, for any irreducible divisor $D$ in $\P^d$ there is a unique irreducible divisor $\widetilde D$ in $Y$, the strict transform of $D$, that maps onto $D$, and the map from $\widetilde D$ to $D$ is birational.

Let $x_0, \ldots, x_d$ be the standard homogeneous coordinates on $\P^d$, and let $H$ be the coordinate hyperplane cut out by $x_0$, with $\widetilde H$ its strict transform in $Y$.  For $1 \leq i \leq n$, let $\chi_i$ be the character of $T$, viewed as a rational function on $X(\Delta_\M)$, corresponding to the unique lattice point $u_i \in \Hom(N,\Z)$ such that
\[
\<u_i, e_j\> = \left\{ \begin{array}{ll} -1 & \mbox{ if } i = 0, \\
								1 & \mbox{ if } i = j, \\
								0 & \mbox{ otherwise.}
\end{array} \right.									
\]
The character $\chi_i$ is invertible on the dense open subset $Y^\circ$ of $Y$, and hence restricts to a nonzero rational function $\chi_i|_Y$ on $Y$.

\begin{lemma} \label{simple pole}
The rational function $\chi_i|_Y$ has a simple pole along $\widetilde H$, and all its other poles are contracted by $p$, for $1 \leq i \leq n$.
\end{lemma}

\begin{proof}
On $X(\Delta_\M)$, the rational function $\chi_i$ has simple poles along those divisors $D_{\rho_F}$ corresponding to flats $F$ that contain $e_0$ but not $e_i$, and no other poles.  In particular, $p$ maps each pole of $\chi_i$ into $H$, and every $T$-invariant divisor of $X(\Delta_\M)$ that maps onto $H$ is a simple pole of $\chi_i$.  It follows that all poles of $\chi_i|_Y$ other than $\widetilde H$ are contracted by $p$, and the multiplicity of the pole along $\widetilde H$ is independent of $i$.  It will therefore suffice to show that $\chi_1|_Y$ has a simple pole along $\widetilde H$.

Now $\chi_1|_Y$ is the pullback of the rational function $x_1/x_0$ on $\P^d$ under the birational map $p$.  Since $x_1/x_0$ has a simple pole along $H$, its pullback has a simple pole along the strict transform $\widetilde H$, and the lemma follows.
\end{proof}

\begin{proof}[Proof of Proposition~\ref{linear subspace}]
To prove the proposition, we must show that the rational function $\chi_i|_Y$ is a linear combination of $1, \chi_1|_Y, \ldots, \chi_d|_Y$, for $i > d$.  Since $p: Y \rightarrow \P^d$ is proper and birational, $\chi_i|_Y$ induces a rational function $f_i$ on $\P^d$ whose poles are the push forward of the poles of $\chi_i|_Y$.  By Lemma~\ref{simple pole}, $f_i$ has a simple pole along $H$ and no other poles, and therefore may be expressed as a linear combination 
\[
f_i = a_0 + a_1 x_1/x_0 + \cdots + a_d x_d/x_0.
\]  
of the $T$-eigensections of $\mathcal{O}(H)$.  Since $\chi_j|_Y$ is the pull back of $x_j/x_0$ for $1 \leq j \leq d$, pulling back the expression above to $Y$ gives $\chi_i|_Y$ as a linear combination of $1, \chi_1|_Y, \ldots, \chi_d|_Y$, as required.
\end{proof}

Let $Y^\circ $ be a tropical realization of $\Delta_\M$, and let $Y$ be its closure in $X(\Delta_\M)$.  By Proposition~\ref{linear subspace}, the image of $Y$ in $\P^n$ is a copy of $\P^d$ embedded linearly, so $Y^\circ$ is exactly the complement in $\P^d$ of the arrangement of $n+1$ hyperplanes given by intersecting with the coordinate hyperplanes in $\P^n$.  In other words, when we factor $p$ as a birational morphism followed by a linear projection, as described above, then $Y^\circ$ is the complement of the hyperplane arrangement $H_0, \ldots, H_n$, where $H_i$ is the projection of $p_1(Y) \cap D_i$, where $D_i$ is the $i$th coordinate hyperplane in $\P^n$.

We now show that this hyperplane arrangement realizes the matroid $\M$, in the following standard sense.  In terms of flats, a collection of hyperplanes $H_0, \ldots, H_n$ realizes $\M$ if the codimension of any intersection $H_{i_0} \cap \cdots \cap H_{i_r}$ is equal to the rank of the flat spanned by $e_{i_0}, \ldots, e_{i_r}$.  This is equivalent to the condition that the intersection of any $d+1$ projective hyperplanes $H_{j_0} \cap \cdots \cap H_{j_d}$ is empty if and only if $\{e_{j_0}, \ldots, e_{j_d} \}$ is a basis for $\M$.  See \cite{Oxley92} for this and other standard facts from matroid theory.

\begin{proposition} \label{realization of M}
Let $H_i$ be the hyperplane given by intersecting $\P^d$ with the $i$th coordinate hyperplane in $\P^n$.  Then $\{H_0, \ldots, H_n\}$ realizes $\M$.
\end{proposition}

\begin{proof}
Let $D_i$ be the $i$th coordinate hyperplane in $\P^n$, so $H_i$ is the intersection of $D_i$ with $\pi(Y) \cong \P^d$.  For any $I \subset \{0, \ldots, n \}$, the preimage in $X(\Delta_\M)$ of the coordinate linear subspace $\bigcap_{i \in I} D_i$ is the union of the $T$-invariant divisors $D_{\rho_F}$ for proper flats $F$ that contain $I$.  Therefore, if $\{e_{i_0}, \ldots, e_{i_d}\}$ is a basis for $\M$ then the preimage of $D_{i_0} \cap \cdots \cap D_{i_n}$ in $X(\Delta_\M)$ is empty.  It follows that the preimage of $H_{i_0} \cap \cdots \cap H_{i_d}$ in $Y$ is empty, and since $\pi:Y \rightarrow \P^d$ is surjective, it follows that $H_{i_0} \cap \cdots \cap H_{i_d}$ is empty as well.

Conversely, if $\{e_{i_0}, \ldots, e_{i_d}\}$ is not a basis for $\M$, then it is a subset of a proper flat $F$.  Now $Y$ meets the divisor $D_{\rho_F}$, and any point in $Y \cap D_{\rho_F}$ projects into $H_{i_0} \cap \cdots \cap H_{i_d}$, so this intersection is nonempty, as required.
\end{proof}

The torus acts naturally on the set of tropical realizations of any fan; if $Y^\circ$ is a tropical realization of $(\Delta,w)$ then the translation $t \cdot Y^\circ$ is also a tropical realization of $(\Delta,w)$, for $t$ in $T$.  We say that two tropical realizations of $\Delta_\M$ are isomorphic if they differ by translation by an element of $T$ and that two realizations of $\M$ are isomorphic if they differ by an automorphism of $\P^d$.
  
\begin{corollary} \label{surjective on points}
The map taking a tropical realization $Y^\circ$ of $\Delta_\M$ to the hyperplane arrangement $H_0, \ldots, H_n$ in $\P^d$ gives a natural bijection between the isomorphism classes of tropical realizations of $\Delta_\M$ over a field $k$ and the isomorphism classes of realizations of $\M$ over $k$.
\end{corollary}

\begin{proof}
First, we show that the map from isomorphism classes of tropical realizations of $\Delta_\M$ to realizations of $\M$ is injective.  Suppose $Y_1^\circ$ and $Y_2^\circ$ are tropical realizations of $\Delta_\M$.  Since there is a unique automorphism of $\P^d$ taking $H_i$ to the $i$th coordinate hyperplane, for $0 \leq i \leq d$, the projected images $p(Y_1^\circ)$ and $p(Y_2^\circ)$ in $\P^d$ must be equal.  Then the two embeddings of this variety in $T$ are determined by the push forwards under $p$ of the rational functions $\chi_i|_{Y_j^\circ}$, for $0 \leq i \leq n$, which are rational functions on $\P^d$ whose divisors of zeros and poles are independent of all choices, as in the proof of Lemma~\ref{simple pole}.  Therefore $\chi_i|_{Y_1^\circ}$ is a nonzero scalar multiple of $\chi_i|_{Y_2^\circ}$.  It follows that $Y_1^\circ$ is the translation of $Y_2^\circ$ by an element of $T$, so the map from isomorphism classes of tropical realizations of $\Delta_\M$ to isomorphism classes of realizations of $\M$ is injective.  It remains to show that every realization of $\M$ comes from a tropical realization of $\Delta_\M$ in this way.

Given a hyperplane arrangement realizing $\M$, there is a natural embedding of the complement of this arrangement in $T$, well-defined up to translation, given by the space of global invertible functions.  The tropicalization of the image is exactly $\Delta_\M$, by the geometric tropicalization theory in Section~2 of \cite{HackingKeelTevelev09}.
\end{proof}

\noindent In particular, there is a natural surjection from the set of $k$-points of $\mathbf R_{\Delta_\M}$ to the set of isomorphism classes of realizations of $\M$ over $k$.  In the following concluding section we show that this surjection is induced by a smooth and surjective morphism of schemes, and deduce that tropical realization spaces satisfy Murphy's Law.

\section{Murphy's Law for tropical realization spaces}

Recall that a hyperplane in $\P^d$ is cut out by a nonzero section of the line bundle $\mathcal{O}(1)$ and a family of hyperplanes over a scheme $S$ is cut out by a section of $\mathcal{O}(1)$ on the relative projective space $\P^d_S$ that is nonzero on the fiber over each point in $S$.  The functor taking a scheme $S$ to the set of such sections is represented by a scheme $\Gamma$ that is isomorphic to $\A^{d+1} \smallsetminus 0$, and there is a locally closed subscheme $\Gamma_\M$ of $\Gamma^{n+1}$ parametrizing tuples of sections $(s_0, \ldots, s_n)$ cutting out families of realizations of the matroid $\M$.  These are exactly the tuples such that $\{ s_{i_0}, \ldots, s_{i_d} \}$ generates $\mathcal{O}(1)$ if and only if $\{e_{i_0}, \ldots, e_{i_d}\}$ is a basis for $\M$.

Note that $\PGL_{d+1}$ acts on $\Gamma_\M$ by automorphisms of $\P^d$ and $\G_m^{n+1}$ acts by coordinatewise scaling of the sections $s_i$.  These actions commute, and the product $\PGL_{d+1} \times \G_m^{n+1}$ acts freely on $\Gamma_\M$.  Furthermore, two tuples of sections cut out the same family of hyperplane arrangements if and only if they differ by an element of $\G_m^{n+1}$ and two families of hyperplane arrangements are isomorphic, by definition, if they differ by an element of $\PGL_{d+1}$.  The quotient
\[
\mathbf R_\M = \Gamma_\M / (\PGL_{d+1} \times \G_m^{n+1})
\]
is the usual \emph{realization space} of the matroid $\M$ and represents the functor taking a scheme $S$ to the set of isomorphism classes of families of hyperplane arrangements realizing $\M$ over $S$.  Theorem~\ref{smooth and surjective} says that $\mathbf R_{\Delta_\M}$ is naturally a torus torsor over $\mathbf R_\M$.  To prove this, we will describe a quotient of $T$ that acts freely on $\mathbf R_{\Delta_\M}$, with quotient $\mathbf R_\M$.

Recall that a \emph{circuit} in a matroid is a minimal set that is not contained in a basis, and a matroid is \emph{connected} if any two elements are contained in a circuit.  The matroid $\M$ has a unique decomposition into connected components 
\[
\M = F_0 \sqcup \cdots \sqcup F_c,
\]
which are disjoint flats such that every circuit of $\M$ is contained in some $F_i$, and the restriction of $\M$ to each $F_i$ is connected.  Let $N_\M$ be the sublattice of $\M$ generated by $v_{F_0}, \ldots, v_{F_c}$.  Note that this sublattice has rank exactly $c$ and is saturated in $N$; the sum $v_{F_0} + \cdots + v_{F_c}$ is zero.  Let $T_\M$ be the subtorus of $T$ whose lattice of one-parameter subgroups is $N_\M$.

\begin{proposition}  \label{subtorus}
The subtorus $T_\M$ acts trivially on the tropical realization space $\mathbf R_{\Delta_\M}$, and the quotient $T/T_\M$ acts freely.
\end{proposition}

To prove the proposition, we use the following lemma characterizing subvarieties of $T$ that are invariant under a subtorus.  Let $T'$ be a subtorus of $T$ whose lattice of one-parameter subgroups is $N'$.

\begin{lemma} \label{translation}
Let $Y^\circ$ be a subvariety of $T$.  Then $Y^\circ$ is invariant under the action of $T'$ if and only if $\Trop(Y^\circ)$ is invariant under translation by $N'_\RR$.
\end{lemma}

\begin{proof}
Suppose $Y^\circ$ is invariant under the action of $T'$.  Then base change to an extension field $K$ with a valuation that surjects onto $\RR$, such as the field of generalized power series $k((t^\RR))$.  Since tropicalization commutes with base change, by Proposition~6.1 of \cite{analytification}, $\Trop(Y^\circ)$ is the image of $Y^\circ(K)$ under the valuation.  Then the action of $T'(K)$ on $Y(K)$ induces an action of $N'_\RR$ on $\Trop(Y^\circ)$ by translations, as required.

Conversely, suppose $\Trop(Y^\circ)$ is invariant under translation by $N'_\RR$.  Then the image of $\Trop(Y^\circ)$ in $N_\RR/N'_\RR$ has dimension $\dim Y^\circ - \dim N'_\RR$.  This image is exactly the tropicalization of the closure of the image of $Y^\circ$ in $T/T'$.  Therefore, the image of $Y^\circ$ in $T/T'$ has dimension $\dim Y^\circ - \dim T'$, and hence $Y^\circ$ is invariant under the action of $T'$.
\end{proof}

\begin{proof}[Proof of Proposition~\ref{subtorus}]
We begin by showing that the matroid fan $\Delta_\M$ is preserved under translation by $(N_\M)_\RR$.  The support of the matroid fan $|\Delta_\M|$ is the image in $N_\RR$ of the set of vectors $(v_0, \ldots, v_n)$ in $\RR^{n+1}$ such that, for every circuit $\{e_{i_0}, \ldots, e_{i_k}\}$ of $\M$, the minimum of the set of coordinates $\{v_{i_0}, \ldots, v_{i_k}\}$ occurs at least twice \cite{FeichtnerSturmfels05}.  Now, let $F$ be a connected component of $\M$.  Since every circuit of $\M$ is either contained in $F$ or disjoint from $F$, it follows that $|\Delta_\M|$ is invariant under translation by $v_F$.  Since this holds for each connected component of $\M$, $|\Delta_\M|$ is invariant under translation by $(N_\M)_\RR$.

Next we show that $(N_\M)_\RR$ is the intersection of the affine spans of the maximal cones in $\Delta_\M$, and hence every vector in $N_\RR$ that acts on $\Delta_\M$ by translations is in $(N_\M)_\RR$.  Suppose the image $v$ in $N_\RR$ of $(v_0, \ldots, v_n) \in \RR^{n+1}$ is in the span of every maximal cone of $\Delta_\M$.  To show that $v$ is in $(N_\M)_\RR$, we must show that $v_j$ is equal to $v_k$ whenever $e_j$ and $e_k$ are in the same connected component of $\M$.  Now, if $e_j$ and $e_k$ are in the same connected component then they are contained in a circuit $\{e_{i_1}, \ldots, e_{i_s}, e_j, e_k \}$.  Then there is a full flag of proper flats
\[
\F = F_1 \subset \cdots \subset F_{d}
\]
such that $F_r$ is the span of $\{e_{i_1}, \ldots, e_{i_r} \}$ for $1 \leq r \leq s$, and $F_{s+1}$ is the span of the full circuit.  The maximal cone $\sigma_\F$ is spanned by $\{ v_{F_1}, \ldots, v_{F_d} \}$, and since each $F_i$ contains either none or both of $e_j$ and $e_k$, and $v$ is in the span of the vectors $v_{F_i}$, it follows that the coordinates $v_j$ and $v_k$ are equal, as required.

By Lemma~\ref{translation}, since $(N_\M)_\RR$ acts on $|\Delta_\M|$ by translations, the subtorus $T_\M$ preserves each tropical realization of $\Delta_\M$, so the action of $T_\M$ on the tropical realization space $\mathbf R_{\Delta_\M}$ is trivial.  It remains to show that $T/T_\M$ acts freely on $\mathbf R_{\Delta_\M}$.

By Fulton's criterion for intersection multiplicity one, the closure of each tropical realization of $\Delta_\M$ meets the orbit $O_\sigma$ corresponding to a maximal cone transversally in a single point.  Therefore, to show that $T/T_\M$ acts freely on $\mathbf R_{\Delta_\M}$ it will suffice to show that $T/T_\M$ acts freely on the union of these orbits $O_\sigma$.  Since the pointwise stabilizer of $O_\sigma$ is exactly the subtorus of $T$ corresponding to the span of $\sigma$, and the intersection of these subtori is $T_\M$, it follows that $T/T_\M$ acts freely on $R_{\Delta_\M}$, as required.
\end{proof}

\noindent We will now prove the following precise version of Theorem~\ref{smooth and surjective}.

\begin{theorem}
The matroid realization space $\mathbf R_\M$ is naturally isomorphic to the quotient of the tropical realization space $\mathbf R_{\Delta_\M}$ by the free action of $T/T_\M$.
\end{theorem}

\begin{proof}

Recall that $\M$ is a fixed matroid on a finite set $E = \{e_0, \ldots, e_n\}$ whose elements are numbered so that $e_0, \ldots, e_d$ is a basis for $\M$.

Let $D_i$ be the $i$th coordinate hyperplane in $\P^n$.  We have defined a natural morphism $p : X(\Delta) \rightarrow \P^d$ that factors as a birational morphism $p_1$ to $\P^n$ followed by the linear projection $p_2$ from $D_0 \cap \cdots \cap D_d$.  If $Y$ is the closure in $X(\Delta_\M)$ of a tropical realization of $\Delta_\M$ then $p_1(Y)$ is a $d$-dimensional linear subspace of $\P^n$ that is disjoint from $D_0 \cap \cdots \cap D_d$, and $p_2$ maps $p_1(Y) \cap D_i$ to a hyperplane $H_i$ in $\P^d$.  In Section~\ref{matroids} we showed that the hyperplane arrangement $H_0, \ldots, H_n$ realizes $\M$.  Now, suppose that $Y_S$ is a family of tropical realizations of $\Delta_\M$ over a scheme $S$.  By the Hilbert polynomial criterion for flatness,
\[
\H_i = p_2(p_1(Y_S) \cap D_i)
\]
is a flat family of hyperplanes.  Let $s_i$ be a section of $\mathcal O(1)$ on $\P^d_S$ cutting out $\H_i$.  By Nakayama's Lemma, $s_{i_0}, \ldots, s_{i_d}$ generate $\mathcal O(1)$ if and only if they generate the fiber at every point.  Therefore it follows from Proposition~\ref{realization of M} that $\H_0, \ldots, \H_n$ is a family of realizations of $\M$, and hence is the pullback under a unique morphism from $S$ to $\mathbf R_\M$.  In particular, the universal family of tropical realizations over $\mathbf R_{\Delta_\M}$ determines a natural morphism $\phi: \mathbf R_{\Delta_\M} \rightarrow \mathbf R_\M$.  We will show that $\phi$ is the quotient morphism for the free action of $T/T_\M$ on $\mathbf R_{\Delta_\M}$ given by Proposition~\ref{subtorus}; to prove this, we first show that the quotient map from $\Gamma_\M$ to $\mathbf R_\M$ factors through $\phi$.

Recall that the closure $Y$ in $X(\Delta_\M)$ of any tropical realization of the matroid fan $\Delta_\M$ over a field is the strict transform of a $d$-dimensional linear subspace of $\P^n$ under a sequence of iterated blowups of strict transforms of coordinate linear subspaces.  We construct a family of tropical realizations of $\Delta_\M$ from a tuple of sections of $\mathcal O(1)$ cutting out a realization of $\M$ by a similar sequence of blowups, giving a natural morphism $\psi: \Gamma_\M \rightarrow \mathbf R_{\Delta_\M}$, as follows.

Let $S$ be a scheme, and let $s_0, \ldots, s_n$ be sections of $\mathcal O(1)$ on $\P^d_S$ cutting out a family of realizations of $\M$.  For each flat $F$ of $\M$, let $\mathcal V_F$ be the family of linear subspaces
\[
\mathcal V_F = \bigcap_{e_i \in F} \H_i
\]
obtained by intersecting the corresponding families of hyperplanes.  We construct a scheme $Y_S$ from a sequence of blowups of $\P^d_S$
\[
Y_S = Y_d \rightarrow Y_{d-1} \rightarrow \cdots \rightarrow Y_1 \rightarrow Y_0 \cong \P^d_S,
\]
where the map from $Y_{j+1}$ to $Y_j$ is the blowup along the strict transform of the flats $\mathcal V_F$ of relative dimension $j$ over $S$.  By construction, $Y_S$ is flat and proper over $S$, and we now construct an embedding of $Y_S$ in $X(\Delta_\M)$, making $Y_S$ a family of tropical realizations of $\Delta_\M$.

By \cite{Cox95b}, to give a morphism from $Y_S$ to $X(\Delta_\M)$ we must give a line bundle $L_F$ for every ray $\rho_F$ in $\Delta_\M$ together with a collection of compatible isomorphisms $\varphi_i: \otimes_F L_F^{\<u, v_F\>} \xrightarrow{\sim} \mathcal O_{Y_S}$ for $u \in \Hom(N, \Z)$.  Since the lattice points $u_i$ considered in Section~\ref{matroids}, for $1 \leq i \leq n$, characterized by the pairings
\[
\<u_i, e_j\> = \left\{ \begin{array}{ll} -1 & \mbox{ if } i = 0, \\
								1 & \mbox{ if } i = j, \\
								0 & \mbox{ otherwise,}
\end{array} \right.									
\]
form a basis for $\Hom(N,\Z)$, it will suffice to give the isomorphisms $\varphi_{u_i}$.  Then all other $\varphi_u$ are given uniquely as a composition of the $\varphi_{u_i}$, and the compatibility condition is automatically satisfied.  Now, the line bundle $\otimes_F L_F^{\<u_i, v_F\>}$ is exactly $\mathcal O(\H_i) \otimes \mathcal O(\H_0)^\vee$, so there are natural isomorphisms given by the sections $s_0/s_i$ of $\Hom(\mathcal O(\H_i) \otimes \mathcal O(\H_0)^\vee, \mathcal O_{Y_S})$.  From the case where $S$ is $\Spec k$, treated in Section~\ref{matroids}, it follows that the resulting map to $X(\Delta_\M) \times S$ is an embedding and that the image is a family of tropical realizations of $\Delta_\M$ over $S$.  By the universal property of the tropical realizations space, this family is pulled back under a unique morphism from $S$ to $\mathbf R_{\Delta_\M}$.  In particular, the universal family of tuples of sections over $\Gamma_\M$ determines a natural morphism $\psi$ from $\Gamma_\M$ to $\mathbf R_{\Delta_\M}$, taking a tuple of sections $(s_0, \ldots, s_n)$ to the tropical realization of $\Delta_\M$ given by embedding the complement of the vanishing loci of the $s_i$ in the dense torus $T$ in $\P^n$ by homogeneous coordinates $[s_0: \cdots :s_n]$.  The $\phi \circ \psi$ takes a tuple of sections $s_0, \ldots, s_n$ to the realization of $\M$ cut out by the $s_i$, and is the quotient map from $\Gamma_\M$ to $\mathbf R_\M$, as claimed.

By construction, $\psi$ is invariant under the action of $PGL_{d+1}$ and the diagonal subtorus $\G_m$ in $\G_m^{n+1}$, and descends to a $T$-equivariant map
\[
\overline \psi: \Gamma_\M / (\PGL_{d+1} \times \G_m) \rightarrow \mathbf R_{\Delta_\M}.
\]
The composition $\phi \circ \overline \psi$ is then the quotient map for a free $T$-action.  Since the subtorus $T_\M$ acts trivially on $\mathbf R_{\Delta_\M}$, and the quotient $T/ T_\M$ acts freely, $\overline \psi$ must be the quotient by the action of $T_\M$ and $\phi$ is the quotient by the action of $T/T_\M$, making $\mathbf R_{\Delta_\M}$ a $T/T_\M$ torsor over $\mathbf R_\M$, as required.
\end{proof}

We conclude by using Theorem~\ref{smooth and surjective} to show that tropical realization spaces satisfy Murphy's Law.

\begin{proof}[Proof of Theorem~\ref{Murphy}]
Let $(Y,y)$ be a singularity of finite type over $\Spec \Z$.  By the scheme theoretic version of Mn\"ev's Universality Theorem \cite{Mnev85}, as presented by Lafforgue in Section~1.8 of \cite{Lafforgue03}, the singularity type of $(Y,y)$ occurs at some point of the realization space of a rank three matroid $\M$.  Since the tropical realization space $\mathbf R_{\Delta_\M}$ is smooth and surjective over $\mathbf R_\M$, by Theorem~\ref{smooth and surjective}, it follows that this singularity type also appears on $\mathbf R_{\Delta_\M}$.  It remains to show that the universal family is smooth over $\mathbf R_{\Delta_\M}$, and since the universal family is flat it is enough to show that the fibers are smooth.  Now the constructions of Section~\ref{matroids} show that each fiber is isomorphic to the blowup of $\P^2$ at finitely many distinct points, and is therefore smooth, as required.
\end{proof}

\bibliography{math}
\bibliographystyle{amsalpha}

\end{document}